\documentclass{article}

\usepackage{graphicx}
\usepackage{amsmath}
\usepackage{a4wide}
\usepackage{cite}

\newtheorem{example}{Example}

\begin{document}

\title{Approximated Analytical Solution to\\
an Ebola Optimal Control Problem\thanks{This is a preprint 
of a paper whose final and definite form is in \emph{International 
Journal for Computational Methods in Engineering Science and Mechanics},
ISSN 1550-2287 (Print), 1550-2295 (Online). Paper Submitted 14-Jul-2015; 
Revised 29-Oct-2015; Accepted for publication 09-Dec-2015.}}


\author{Doracelly Hincapi\'{e}-Palacio$^1$\\
{\small {\tt doracelly.hincapie@udea.edu.co}}
\and Juan Ospina$^2$\\
{\small {\tt jospina@eafit.edu.co}}
\and Delfim F. M. Torres$^3$\\
{\small {\tt delfim@ua.pt}}}


\date{$^1$Group of Epidemiology, National School of Public Health\\
University of Antioquia, Medellin, Colombia\\[0.3cm]
$^2$Logic and Computation Group, School of Sciences\\
Universidad EAFIT, Medellin, Colombia\\[0.3cm]
$^3$\text{Center for Research \& Development in Mathematics and Applications (CIDMA)}\\
Department of Mathematics, University of Aveiro, 3810-193 Aveiro, Portugal}


\maketitle


\begin{abstract}
An analytical expression for the optimal control
of an Ebola problem is obtained. The analytical
solution is found as a first-order approximation
to the Pontryagin Maximum Principle via the Euler--Lagrange equation.
An implementation of the method is given using the computer algebra system
\textsf{Maple}. Our analytical solutions confirm the results recently reported
in the literature using numerical methods.

\medskip

\noindent {\bf Keywords:} Optimal Control; Euler--Lagrange equation;
Computer Algebra; Ebola; Approximated analytical expressions.

\smallskip

\noindent {\bf Mathematics Subject Classification 2010:} 49-04; 49K15; 92D30.
\end{abstract}


\section{Introduction}
\label{intro}

The largest outbreak of Ebola virus ever recorded has been ongoing
since was first confirmed in March, 2014.
Ebola is a fatal disease that has claimed 7 000 lives
by the end of 2014 in just Guinea, Liberia and Sierra Lione.
While the Ebola outbreak has slowed down across West Africa by June 2015,
every new infection continues to threaten millions of lives
and bringing fear to the world. With more than 24 000 cases
and almost 10 000 fatalities, this outbreak is already one
of the biggest public health crises of the XXI century.
Overcoming Ebola is a complex emergency, challenging
not only governments and international aid organisations
but also computational and life scientists and applied mathematicians
\cite{Atangana,Lewnard,Marzi,MyID:321,MyID:322,WangZhong}.

In a recent work by Rachah and Torres, an optimal control
problem of the 2014 Ebola outbreak in West Africa was
posed and numerically solved through \textsf{Matlab}
and the \textsf{ACADO} toolkit \cite{MyID:321}.
See also \cite{MyID:322} for a different model and other
\textsf{Matlab} numerical simulations. In contrast,
here we address the problem by analytical methods.
The results confirm the previous numerical results,
but now with a theoretical/analytical foundation.
The new method is simple but envolves lengthy calculations.
For this reason, a computer algebra package with the proposed method
is developed in \textsf{Maple}.

The text is organized as follows. In Section~\ref{prob}
the optimal control problem is formulated. Our method is
explained in Section~\ref{meth} and illustrated with an example.
Then, in Section~\ref{resu}, we apply it to the Ebola optimal
control problem. We end with Section~\ref{conclu} of
conclusions, while Appendix~\ref{apen} provides the
developed \textsf{Maple} code.


\section{The Problem}
\label{prob}

The Ebola problem of optimal control proposed in \cite{MyID:321}
consists to determine the control function $u(\cdot)$ in such
a way the objective functional given by
\begin{equation}
\label{eq:E4}
\mathcal{J}\left( u \right) =
\int_{0}^{T} \!\left[I\left(t\right)
+\frac{1}{2}A u \left( t \right)^{2}\right]{dt}
\end{equation}
is minimized, where $A$ is a fixed nonnegative constant,
when subject to the dynamic equations
\begin{gather}
{\frac {d}{dt}}S \left( t \right) =-\beta\,S \left( t \right) I
 \left( t \right) -u \left( t \right) S \left( t \right) \label{eq:E1} \\
{\frac {d}{dt}}I\left( t \right) =\beta\,S \left( t \right) I \left(
t \right) -\mu\,I \left( t \right) \label{eq:E2}\\
{\frac {d}{dt}}R \left( t \right) =\mu\,I \left( t \right)
+u \left( t\right) S \left( t \right) \label{eq:E3}
\end{gather}
for all
\begin{equation}
\label{eq:E7}
t \in [0,T],
\end{equation}
the given initial conditions
\begin{equation}
\label{eq:E5}
S(0) \geq 0, \quad  I(0) \geq 0, \quad R(0) \geq 0,
\end{equation}
and where the control values are bounded in the interval $[0, 0.9]$, that is,
\begin{equation}
\label{eq:E6}
0 \leq u(t) \leq  0.9.
\end{equation}
Here $T$ is the duration of the application of the control (duration of the
vaccination program). The constant $0.9$ is a control value that
is able to eliminate the Ebola transmission according with $R_{0} <1$,
where $R_{0}$ is the basic reproduction number for the system
\eqref{eq:E1}--\eqref{eq:E3}, and control $u$ is considered constant
along all time. This control $u(t) \equiv 0.9$ is however not optimal.
For this reason, we search for an optimal value of $u(t)$, $t \in [0, T]$,
subject to the constraint given by \eqref{eq:E6}. Note that the control $u(t)$
represents the vaccination rate at time $t$. Being the fraction of susceptible
individuals vaccinated per unit of time, the value $0.9$ means that,
at maximum, 90\% of susceptible are vaccinated. In other words, what we assume
here is that the fraction of individuals who are not vaccinated takes
at least the value of 10\%. This is in agreement with general experience
in vaccination, where it is well recognized the impossibility
to vaccinate all population. For more details on the
description of the mathematical model and the meaning of the parameters,
we refer the reader to the work of Rachah and Torres \cite{MyID:321,MyID:322}.
In particular, see the scheme of the susceptible-infected-recovered model
with vaccination found in Section~4.2 of \cite{MyID:321}
and the optimal control problem in Section~5 of \cite{MyID:321},
where the following parameters, initial conditions, and time horizon
are considered: infection rate $\beta = 0.2$; recovery rate $\mu = 0.1$;
at the beginning 95\% of population is susceptible and 5\% is already infected,
that is, $S(0) = 0.95$, $I(0) = 0.05$ and $R(0) = 0$; and $T = 100$ days.
Differently from previous works \cite{MyID:321,MyID:322},
which are exclusively based on numerical methods,
we address the optimal control problem \eqref{eq:E4}--\eqref{eq:E6}
by using an approximated analytic method. For that we make use
of the computer algebra system \textsf{Maple}.


\section{The Method}
\label{meth}

In this section the approximated analytical method that is used
in Section~\ref{resu} to solve the optimal control problem
\eqref{eq:E4}--\eqref{eq:E6} is explained and illustrated
with an example. The idea is to use the
classical calculus of variations, specifically
the Euler--Lagrange equation, which is its main tool.
The Euler--Lagrange equation is used with the aim to obtain
a first-order approximation to the Pontryagin Maximum
Principle. Typically, the Pontryagin Maximum Principle is harder
to solve analytically than the Euler--Lagrange equation. In contrast,
the Euler--Lagrange equations can be easily solved analytically
in many interesting cases. In our work we perform an analytical experiment
consisting to solve analytically the optimal control problem
\eqref{eq:E4}--\eqref{eq:E6}, which was previously solved numerically
in \cite{MyID:321}. As we shall see in Section~\ref{resu},
our approach turns out to be a good one.

Let us start with the dynamical control system
\begin{equation}
\label{eq:E8}
{\frac {d}{dt}}y \left( t \right)
=F\left( y \left( t \right), u\left(t\right)\right),
\end{equation}
where $y(\cdot)$ is the state vector that must be controlled
and $u(\cdot)$ is the control that must be applied to the system
in order to minimize the functional
\begin{equation}
\label{eq:E9}
\mathcal{I}\left( u \right) =\int _{0}^{T} \left[
\left( {\frac {d}{dt}}u \left( t \right)\right)^{2}
+\frac{1}{2} A u\left(t \right)^{2}\right]{d t},
\end{equation}
where $A$ is the  parameter that is determining the cost of the control
and $T$ is the duration of application of the control. Comparing
the objective functionals \eqref{eq:E4} and \eqref{eq:E9}, we are assuming that
$I(t)=(du(t)/dt)^2$. This is a particular case of the more general assumption
\begin{equation}
\label{series:I}
I(t)=a_{1}(t)u(t)+a_{2}(t)(du(t)/dt)^2+a_{3}(t)(du(t)/dt)^4+\cdots.
\end{equation}
From the epidemiological point of view, given that system \eqref{eq:E1}--\eqref{eq:E3}
can be considered as a black box, being the input $u(t)$ and the output $I(t)$,
it is possible to think that $I(t)$ is approximately given by a series of the
form \eqref{series:I}. Given that $I(t)$ is always positive, we use even powers
of $du(t)/dt$. The simplest assumption is then $I(t)=(du(t)/dt)^2$, which makes
functional \eqref{eq:E9} to take the form of the Lagrangian for the classical
harmonic oscillator. It is possible to use other forms for $I(t)$ as a function
of $u(t)$ and its derivatives. For our purposes, the simplest expression
$I(t)=(du(t)/dt)^2$ is enough. The Euler--Lagrange equation
(see, e.g., \cite{MR3156230}) associated with \eqref{eq:E9} is
\begin{equation}
\label{eq:E10}
Au \left( t \right) -2\,{\frac {d^{2}}{d{t}^{2}}}u \left( t \right) =0
\end{equation}
and the solution of \eqref{eq:E10} with the conditions
\begin{equation}
\left\{ u \left( 0 \right) =U_{{0}},u \left( \infty  \right) =0 \right\}
\end{equation}
is given by
\begin{equation}
\label{eq:E12}
u \left( t \right) =U_{{0}}{{\rm e}^{-\frac{1}{2}\sqrt {2}\sqrt {A}t}}.
\end{equation}
In \eqref{eq:E10} we are assuming that $u$
is of class $C^2$: the classical Euler--Lagrange equation is a
second-order differential equation. The exact solution
$u$ is not necessarily $C^2$, but it can always be
approximated by a $C^2$ function. Note that our goal
is to find an approximated analytical solution and not the exact one.
Replacing \eqref{eq:E12} in \eqref{eq:E8}, we obtain that
\begin{equation}
\label{eq:E13}
{\frac {d}{dt}}y \left( t \right)
=F \left( y \left( t \right) ,U_{{0}}{{\rm e}^{-\frac{1}{2}
\sqrt{2}\sqrt{A}t}} \right).
\end{equation}
Now we assume that equation \eqref{eq:E13} can be solved analytically
when subject to the initial condition $y(0) = Y_0$. Then, formally,
it is possible to write that
\begin{equation}
\label{eq:E14}
y \left( t \right) =G \left( t,A,U_{{0}},Y_{{0}} \right).
\end{equation}
To determine $U_0$, we minimize the following functional:
\begin{equation}
\label{eq:E15}
\mathcal{K}\left( u \right)
=\int_{0}^{T}\!\left[y \left( t \right)
+\frac{1}{2}A \left( u \left( t \right)\right)^{2}\right]{dt}.
\end{equation}
Replacing \eqref{eq:E12} and \eqref{eq:E14} in \eqref{eq:E15}, we obtain that
\begin{equation}
\label{eq:E16}
K \left( U_{{0}} \right)
=\int _{0}^{T}\!G \left( t,A,U_{{0}},Y_{{0}} \right) {dt}
-\frac{1}{4}\sqrt {A}\sqrt {2}{U_{{0}}}^{2} \left( -1+{
{\rm e}^{-\sqrt {2}\sqrt {A}T}} \right).
\end{equation}
Taking the derivative of \eqref{eq:E16} with respect to $U_0$
and equating the result to zero, we have
\begin{equation}
\label{eq:E17}
\int_{0}^{T}\!{\frac {\partial }{\partial U_{{0}}}}
G\left( t,A,U_{{0}},Y_{{0}} \right) {dt}
-\frac{1}{2}\sqrt {A}\sqrt {2}U_{{0}} \left( -1
+{{\rm e}^{-\sqrt {2}\sqrt {A}T}} \right) =0.
\end{equation}
The parameter $A$ is determined according to
\begin{equation}
\label{eq:E18}
U_{{0}}{{\rm e}^{-\frac{1}{4}\sqrt {2}\sqrt {A}T}}={\frac {U_{{0}}}{Q}},
\end{equation}
that is, we assume that at the half of the duration of the application of the control,
the intensity of the control is reduced by a factor $Q$ with respect to the
initial intensity. Then the solution of \eqref{eq:E18} is given by
\begin{equation}
\label{eq:E19}
A=8\,{\frac { \left( \ln  \left( Q \right)  \right) ^{2}}{{T}^{2}}}.
\end{equation}
Finally, solving equation \eqref{eq:E17} with respect to $U_0$,
using \eqref{eq:E19} and the numerical values for the other parameters,
the values for $U_0$ and $A$ are obtained and the explicit
form of the control $u(t)$ given by \eqref{eq:E12} is specified.

To illustrate the method that was just explained,
we consider now a simple toy model.

\begin{example}
\label{ex}
Let
\begin{equation}
{\frac {d}{dt}}S \left( t \right) =-\beta\,S \left( t \right) I
 \left( t \right)
\end{equation}
and
\begin{equation}
\label{eq:E21}
{\frac {d}{dt}}I \left( t \right) =\beta\,S \left( t \right) I \left(
t \right) -u \left( t \right) I \left( t \right).
\end{equation}
The problem here is to control the variable $I(t)$ using $u(t)$.
We assume that the control $u(t)$ has the form given by \eqref{eq:E12}.
The expression \eqref{eq:E12} is the solution of the
differential equation \eqref{eq:E10}, which is the Euler--Lagrange equation
for the functional \eqref{eq:E9} with the assumption $I(t)=(du(t)/dt)^2$.
If the more general assumption \eqref{series:I} is used, then the corresponding
Euler--Lagrange equation will be more complex and the explicit solution will
involve special functions, such as Airy, Bessel, Kummer, Whittaker,
and Heun functions. Replacing \eqref{eq:E12} in \eqref{eq:E21}, we obtain that
\begin{equation}
\label{eq:E22}
{\frac {d}{dt}}I\left( t \right) =\beta\,S \left( t \right) I \left(
t \right) -U_{{0}}{{\rm e}^{-\frac{1}{2}\sqrt {2}\sqrt {A}t}}I \left( t\right).
\end{equation}
An approximated analytical solution of equation \eqref{eq:E22} can be
obtained for the early stages of the outbreak when $S(t) \approx S_0$.
With this approximation, \eqref{eq:E22} is reduced to
\begin{equation}
\label{eq:E23}
{\frac {d}{dt}}I\left( t \right) =\beta\,S_{0} I \left(
t \right) -U_{{0}}{{\rm e}^{-\frac{1}{2}\sqrt {2}\sqrt {A}t}}I \left( t
 \right)
\end{equation}
and the explicit solution of \eqref{eq:E23}
with the initial condition $I(0)=i_0$ is given by
\begin{equation}
\label{eq:E24}
I \left( t \right) =i_{{0}}{{\rm e}^{{\frac {-U_{{0}}\sqrt {2}+\beta\,
S_{{0}}t\sqrt {A}+U_{{0}}\sqrt {2}{{\rm e}^{-\frac{1}{2}\sqrt {2}\sqrt {A}t}
}}{\sqrt {A}}}}}.
\end{equation}
For the early stages of the outbreak, equation \eqref{eq:E24} takes the form
\begin{equation}
\label{eq:E25}
I\left( t \right) =i_{{0}} \left( 1+\beta\,S_{{0}}t-tU_{{0}} \right).
\end{equation}
Using \eqref{eq:E15} with $y(t)=I(t)$ and \eqref{eq:E25}, we derive that
\begin{equation}
\label{eq:E26}
K \left( U_{{0}} \right) =i_{{0}}T+\frac{1}{2}i_{{0}}{T}^{2}\beta\,S_{{0}}
-\frac{1}{2}i_{{0}}{T}^{2}U_{{0}}+\frac{1}{4}\sqrt {2}\sqrt {A}{U_{{0}}}^{2}
-\frac{1}{4}\sqrt {2}\sqrt {A}{U_{{0}}}^{2}{{\rm e}^{-\sqrt {2}\sqrt {A}T}}.
\end{equation}
Taking the derivative of \eqref{eq:E26} with respect to $U_0$,
equating the result to zero and solving with respect to $U_0$, we have that
\begin{equation}
\label{eq:E27}
U_{{0}}=-\frac{1}{2}{\frac {i_{{0}}{T}^{2}\sqrt {2}}{\sqrt {A}
\left( -1+{{\rm e}^{-\sqrt {2}\sqrt {A}T}} \right)}}.
\end{equation}
The control $u(t)$ is completely determined by replacing
\eqref{eq:E27} and \eqref{eq:E19} in \eqref{eq:E12}.
All these computations are easily done with the help
of a computer algebra system (see Appendix~\ref{sec:ap:a1}).
\end{example}


\section{Main Results}
\label{resu}

With the aim to apply the method explained in Section~\ref{meth}
to the Ebola problem \eqref{eq:E4}--\eqref{eq:E6}, we assume
that equation \eqref{eq:E1} can be reduced to
\begin{equation}
\label{eq:E28}
{\frac {d}{dt}}S \left( t \right) =-U_{{0}}{{\rm e}^{-\frac{1}{2}\sqrt {2}
\sqrt {A}t}}S \left( t \right)
\end{equation}
at the very early stages of the outbreak. In other words, we assume that
$\beta S(t)I(t) \ll u(t)S(t)$ for $t$ near to zero, that is, at the
beginning of the outbreak the depletion in the number of susceptible
individuals is due to the vaccination, given that the reduction in the number
of susceptible individuals due to infection is depreciated. Then the solution
of \eqref{eq:E28} with initial condition $S(0) = S_0$  is given by
\begin{equation}
\label{eq:E29}
S \left( t \right) =S_{{0}}{{\rm e}^{{\frac {\sqrt {2} U_{{0}}\left(-1
+{{\rm e}^{-\frac{1}{2}\sqrt {2}\sqrt{A}t}} \right)}{\sqrt{A}}}}}.
\end{equation}
At the beginning of the outbreak, \eqref{eq:E29} is reduced to
\begin{equation}
\label{eq:E30}
S \left( t \right) =S_{{0}}-S_{{0}}U_{{0}}t+S_{{0}} \left( \frac{1}{4}
\sqrt{2}U_{{0}}\sqrt {A}+\frac{1}{2}{U_{{0}}}^{2} \right) {t}^{2}.
\end{equation}
Now equation \eqref{eq:E2} with \eqref{eq:E29} takes the form
\begin{equation}
\label{eq:E31}
{\frac {d}{dt}}I \left( t \right) =\beta\,S_{{0}}{{\rm e}^{{\frac{
\sqrt {2}U_{{0}} \left( -1+{{\rm e}^{-\frac{1}{2}\sqrt {2}\sqrt {A}t}}
\right)}{\sqrt {A}}}}}I \left( t \right) -\mu\,I \left( t \right).
\end{equation}
The solution of \eqref{eq:E31} with initial condition
$I(0)=i_0$  is given by
\begin{equation}
\label{eq:E32}
I \left( t \right) =\frac{i_{{0}}{{\rm e}^{ \left( \beta\,
S_{{0}}\sqrt {2}{{\rm e}^{-{\frac{\sqrt {2}U_{{0}}}{\sqrt {A}}}}}{\it Ei}
\left( 1,-{\frac {\sqrt {2}U_{{0}}{{\rm e}^{-\frac{1}{2}
\sqrt{2}\sqrt {A}t}}}{\sqrt {A}}} \right) -\mu\,t
\sqrt {A} \right) {\frac {1}{\sqrt {A}}}}}}{{{\rm e}^{\beta\,S_{{0}}
\sqrt{2}{{\rm e}^{-{\frac {\sqrt {2}U_{{0}}}{\sqrt{A}}}}}{\it Ei} \left( 1,
-{\frac {\sqrt {2}U_{{0}}}{\sqrt {A}}}\right) {\frac {1}{\sqrt {A}}}}}},
\end{equation}
where $Ei(x)$ is the exponential integral function defined by
\begin{equation}
\label{eq:E33}
{\it Ei} \left( x \right)
=-\int _{-x}^{\infty }\!{\frac {{{\rm e}^{-t}}}{t}}{dt}.
\end{equation}
For the early stages of the outbreak, equality \eqref{eq:E32} is reduced to
\begin{equation}
\label{eq:E34}
I\left( t \right) =i_{{0}}+i_{{0}} \left( \beta\,S_{{0}}-\mu \right)
t+C_{{2}}{t}^{2}+\frac{i_0}{12}C_{{3}}{t}^{3}-\frac{i_0}{48}C_{{4}}{t}^{4},
\end{equation}
where
\begin{equation}
\label{eq:E35}
C_{{2}}=i_{{0}} \left( -\frac{1}{2}\beta\,S_{{0}}U_{{0}}
+\frac{1}{2}{\beta}^{2}{S_{{0}}}^{2}-\beta\,S_{{0}}\mu
+\frac{1}{2}{\mu}^{2} \right),
\end{equation}
\begin{equation}
\label{eq:E36}
\begin{split}
C_{{3}} = \beta\,S_{{0}} &U_{{0}}\sqrt {2}\sqrt {A}+2\,\beta\,
S_{{0}}{U_{{0}}}^{2}-6\,{\beta}^{2}{S_{{0}}}^{2}U_{{0}}\\
&+6\,\beta\,S_{{0}}U_{{0}}\mu+2\,{\beta}^{3}{S_{{0}}}^{3}
-6\,{\beta}^{2}{S_{{0}}}^{2}\mu\\
&+ 6\,\beta\,S_{{0}}{\mu}^{2}-2\,{\mu}^{3},
\end{split}
\end{equation}
and
\begin{equation}
\label{eq:E37}
\begin{split}
C_4 =  \beta\, & S_{{0}}U_{{0}}A+3\,\beta\,S_{{0}}{U_{{0}}}^{2}
\sqrt {2}\sqrt {A}+2\,\beta\,S_{{0}}{U_{{0}}}^{3}\\
&  - 4\,{\beta}^{2}{S_{{0}}}^{2}U_{{0}}\sqrt {2}\sqrt {A}
-14\,{\beta}^{2}{S_{{0}}}^{2}{U_{{0}}}^{2}\\
& + 4\,\beta\,S_{{0}}U_{{0}}\mu\,\sqrt {2}\sqrt {A}
+8\,\beta\,S_{{0}}{U_{{0}}}^{2}\mu+12\,{\beta}^{3}{S_{{0}}}^{3}U_{{0}}
-24\,{\beta}^{2}{S_{{0}}}^{2}U_{{0}}\mu\\
& + 12\,\beta\,S_{{0}}U_{{0}}{\mu}^{2}-2\,{\beta}^{4}{S_{{0}}}^{4}
+8\,{\beta}^{3}{S_{{0}}}^{3}\mu-12\,{\beta}^{2}{S_{{0}}}^{2}{\mu}^{2}\\
& + 8\,\beta\,S_{{0}}{\mu}^{3}-2\,{\mu}^{4}.
\end{split}
\end{equation}
Replacing \eqref{eq:E34}--\eqref{eq:E37} and \eqref{eq:E12}
into the functional \eqref{eq:E15}, with $y(t)=I(t)$, we obtain that
\begin{multline}
\label{eq:E38}
K \left( U_{{0}} \right) =  i_{{0}}T+E_{{2}}{T}^{2}+E_{{3}}{T}^{3}
+\frac{1}{48}i_{{0}}E_{{4}}{T}^{4}-{\frac {1}{240}}\,i_{{0}}E_{{5}}{T}^{5}\\
+ \frac{1}{4}{U_{{0}}}^{2}\sqrt {2}\sqrt {A}-\frac{1}{4}\sqrt {2}
\sqrt {A}{U_{{0}}}^{2}{{\rm e}^{-\sqrt {2}\sqrt {A}T}},
\end{multline}
where
\begin{equation}
\label{eq:E39}
E_{{2}}=\frac{1}{2}i_{{0}} \left( \beta\,S_{{0}}-\mu \right),
\end{equation}
\begin{equation}
\label{eq:E40}
E_{{3}}=\frac{1}{3}i_{{0}} \left( -\frac{1}{2}\beta\,S_{{0}}U_{{0}}
+\frac{1}{2}{\beta}^{2}{S_{{0}}}^{2}
-\beta\,S_{{0}}\mu+\frac{1}{2}{\mu}^{2} \right),
\end{equation}
\begin{equation}
\label{eq:E41}
E_4 = \beta\,S_{{0}}U_{{0}}\sqrt {2}\sqrt {A}+2\,\beta\,S_{{0}}{U_{{0}}}^{2}
-6\,{\beta}^{2}{S_{{0}}}^{2}U_{{0}}+6\,\beta\,S_{{0}}U_{{0}}\mu
+2\,{\beta}^{3}{S_{{0}}}^{3}
-6\,{\beta}^{2}{S_{{0}}}^{2}\mu+6\,\beta\,S_{{0}}{\mu}^{2}-2\,{\mu}^{3},
\end{equation}
and
\begin{equation}
\label{eq:E42}
\begin{split}
E_5 = \beta & S_{{0}}U_{{0}}A+3\,\beta\,S_{{0}}{U_{{0}}}^{2}\sqrt {2}\sqrt {A}
+2\,\beta\,S_{{0}}{U_{{0}}}^{3}-4\,{\beta}^{2}{S_{{0}}}^{2}U_{{0}}
\sqrt {2}\sqrt {A}-14\,{\beta}^{2}{S_{{0}}}^{2}{U_{{0}}}^{2}\\
&+ 4\,\beta\,S_{{0}}U_{{0}}\mu\,\sqrt {2}\sqrt {A}+8\,\beta\,
S_{{0}}{U_{{0}}}^{2}\mu+12\,{\beta}^{3}{S_{{0}}}^{3}U_{{0}}
-24\,{\beta}^{2}{S_{{0}}}^{2}U_{{0}}\mu+12\,\beta\,S_{{0}}U_{{0}}{\mu}^{2}\\
&- 2\,{\beta}^{4}{S_{{0}}}^{4}+8\,{\beta}^{3}{S_{{0}}}^{3}\mu
-12\,{\beta}^{2}{S_{{0}}}^{2}{\mu}^{2}+8\,\beta\,S_{{0}}{\mu}^{3}-2\,{\mu}^{4}.
\end{split}
\end{equation}
Taking the derivative of \eqref{eq:E38} with respect to $U_0$,
using \eqref{eq:E39}--\eqref{eq:E42}, and equating the result to zero,
we have that
\begin{equation}
\label{eq:E43}
-\frac{1}{6}i_{{0}}\beta\,S_{{0}}{T}^{3}+F_{{4}}{T}^{4}-{\frac {1}{240}}\,
i_{{0}}F_{{5}}{T}^{5}+\frac{1}{2}\sqrt {2}U_{{0}}\sqrt {A}-\frac{1}{2}\sqrt{2}
\sqrt {A}U_{{0}}{{\rm e}^{-\sqrt {2}\sqrt {A}T}}=0,
\end{equation}
where
\begin{equation}
\label{eq:E44}
F_{{4}}=\frac{1}{48}i_{{0}} \left( \beta\,S_{{0}}\sqrt {2}\sqrt {A}+4\,\beta
\,S_{{0}}U_{{0}}-6\,{\beta}^{2}{S_{{0}}}^{2}+6\,\beta\,S_{{0}}\mu
\right)
\end{equation}
and
\begin{equation}
\label{eq:E45}
\begin{split}
F_5 = \beta & S_{{0}}A+6\,\beta\,S_{{0}}U_{{0}}\sqrt {2}\sqrt {A}
+6\,\beta\,S_{{0}}{U_{{0}}}^{2}-4\,{\beta}^{2}{S_{{0}}}^{2}\sqrt {2}\sqrt{A}
- 28\,{\beta}^{2}{S_{{0}}}^{2}U_{{0}}\\
&+ 4\,\beta\,S_{{0}}\mu\,\sqrt {2}\sqrt {A}+16\,\beta\,S_{{0}}U_{{0}}\mu
+12\,{\beta}^{3}{S_{{0}}}^{3}-24\,{\beta}^{2}{S_{{0}}}^{2}\mu
+12\,\beta\,S_{{0}}{\mu}^{2}.
\end{split}
\end{equation}
Solving \eqref{eq:E43} with respect to $U_0$ and taking into account
\eqref{eq:E44}--\eqref{eq:E45}, we derive that
\begin{equation}
\label{eq:E45A}
U_{{0}}=-{\frac {V-\sqrt {W}}{6 i_{{0}}{T}^{5}\beta\,S_{{0}}}},
\end{equation}
where
\begin{multline}
\label{eq:E45B}
V = -14\,i_{{0}}{T}^{5}{\beta}^{2}{S_{{0}}}^{2}+3\,i_{{0}}{T}^{5}
\beta\,S_{{0}}\sqrt {2}\sqrt {A}+8\,i_{{0}}{T}^{5}\beta\,S_{{0}}\mu \\
- 60\,\sqrt {2}\sqrt {A}+60\,\sqrt {2}\sqrt {A}{{\rm e}^{-\sqrt {2}\sqrt {A} T}}
-10\,i_{{0}}{T}^{4}\beta\,S_{{0}}
\end{multline}
and
\begin{equation}
\label{eq:E45C}
\begin{split}
W &= 960 i_{{0}}{T}^{5}\beta\,S_{{0}}\mu\,\sqrt {2}
\sqrt{A}{{\rm e}^{-\sqrt {2}\sqrt {A}T}}+20\,{i_{{0}}}^{2}{T}^{9}{
\beta}^{2}{S_{{0}}}^{2}\mu-1200\,\sqrt {2}\sqrt {A}{{\rm e}^{
-\sqrt {2}\sqrt {A}T}}i_{{0}}{T}^{4}\beta\,S_{{0}}\\
& + 12\,{i_{{0}}}^{2}{T}^{10}{\beta}^{2}{S_{{0}}}^{2}A
+124\,{i_{{0}}}^{2}{T}^{10}{\beta}^{4}{S_{{0}}}^{4}+100\,
{i_{{0}}}^{2}{T}^{9}{\beta}^{3}{S_{{0}}}^{3}-140
\,{i_{{0}}}^{2}{T}^{8}{\beta}^{2}{S_{{0}}}^{2}+7200\,A\\
& + 1200\,i_{{0}}{T}^{4}\beta\,S_{{0}}\sqrt {2}\sqrt {A}+1680\,
i_{{0}}{T}^{5}{\beta}^{2}{S_{{0}}}^{2}\sqrt {2}\sqrt {A}-960\,
i_{{0}}{T}^{5}\beta\,S_{{0}}\mu\,\sqrt {2}\sqrt {A}\\
&-14400\,A{{\rm e}^{-\sqrt {2}\sqrt {A}T}}
+ 7200\,A \left( {{\rm e}^{-\sqrt {2}\sqrt {A}T}}\right)^{2}
-60\,{i_{{0}}}^{2}{T}^{10}{\beta}^{3}{S_{{0}}}^{3}\sqrt{2}\sqrt {A}
-80\,{i_{{0}}}^{2}{T}^{10}{\beta}^{3}{S_{{0}}}^{3}\mu\\
& -1680\,i_{{0}}{T}^{5}{\beta}^{2}{S_{{0}}}^{2}\sqrt {2}
\sqrt {A}{{\rm e}^{-\sqrt {2}\sqrt {A}T}}
+ 24\,{i_{{0}}}^{2}{T}^{10}{\beta}^{2}{S_{{0}}}^{2}\sqrt {2}\sqrt {A}\mu
-720\,i_{{0}}{T}^{5}\beta\,S_{{0}}A\\
& +720\,i_{{0}}{T}^{5}\beta\, S_{{0}}A{{\rm e}^{-\sqrt {2}\sqrt {A}T}}
-30\,{i_{{0}}}^{2}{T}^{9}{\beta}^{2}{S_{{0}}}^{2}\sqrt {2}\sqrt {A}
-8\,{i_{{0}}}^{2}{T}^{10}{\beta}^{2}{S_{{0}}}^{2}{\mu}^{2}.
\end{split}
\end{equation}
Now we use the following numerical values for the relevant parameters:
\begin{equation}\label{eq:E46}
\left\{ Q=500,T=100,\mu= 0.1,\beta= 0.2,S_{{0}}= 0.95, i_{{0}}= 0.05 \right\}.
\end{equation}
These values are used here for numerical experimentation.
It is, however, possible to consider other values (the concrete values
are not critical for the experiments). We obtain from \eqref{eq:E19} that
\begin{equation}
\label{eq:E47}
A= 0.03089708305.
\end{equation}
Using \eqref{eq:E47}, \eqref{eq:E46} and the expression
for $U_0$ given by \eqref{eq:E45A}--\eqref{eq:E45C}, we obtain that
\begin{equation}
\label{eq:E48}
U_{{0}}= 0.3796479517.
\end{equation}
Replacing \eqref{eq:E47} and \eqref{eq:E48} in \eqref{eq:E12},
we obtain that the optimal control is given by
\begin{equation}
\label{eq:E49}
u \left( t \right) = 0.3796479517\,{{\rm e}^{- 0.1242921619\,t}}.
\end{equation}
Now the system \eqref{eq:E1}--\eqref{eq:E3}
is numerically solved with \eqref{eq:E49} and the initial conditions
\begin{equation}
\label{eq:E50}
\left\{S \left( 0 \right) = 0.95,
 I\left( 0 \right) = 0.05, R \left( 0 \right) =0 \right\}.
\end{equation}
We obtain the curves of Figures~\ref{fig:6}--\ref{fig:9} (for all the details
see the \textsf{Maple} code in Appendix~\ref{sec:ap:a2}).
Our results reproduce the numerical results
of Rachah and Torres \cite{MyID:321} using the simplest assumption
$I(t)=(du(t)/dt)^2$. Note that this assumption is not directly linked
to the numerical results of \cite{MyID:321}: the assumption
$I(t)=(du(t)/dt)^2$ is a particular case of \eqref{series:I}.
In the case we do not have the numerical results in advance,
we can use the general form \eqref{series:I} and experiment
with different terms of such series to get the best possible results.
\begin{figure}
\includegraphics[angle=-90,scale=0.6]{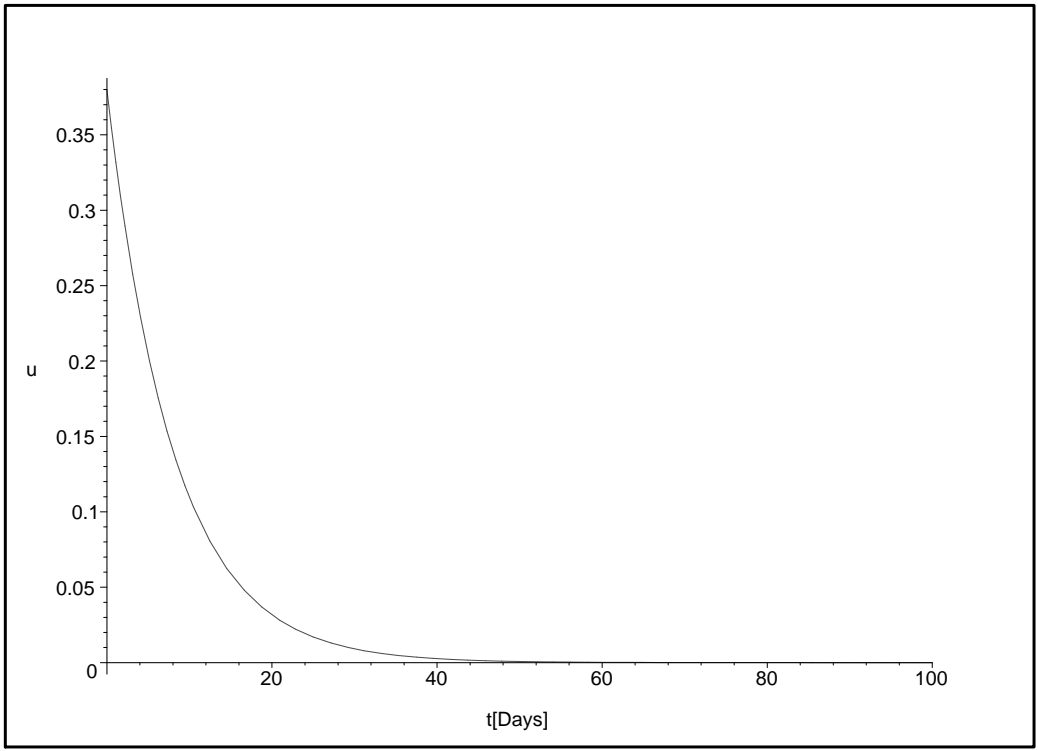}
\caption{Analytical optimal control $u(t)$ \eqref{eq:E49}
for problem \eqref{eq:E4}--\eqref{eq:E6}.}
\label{fig:6}
\end{figure}
\begin{figure}
\includegraphics[angle=-90,scale=0.6]{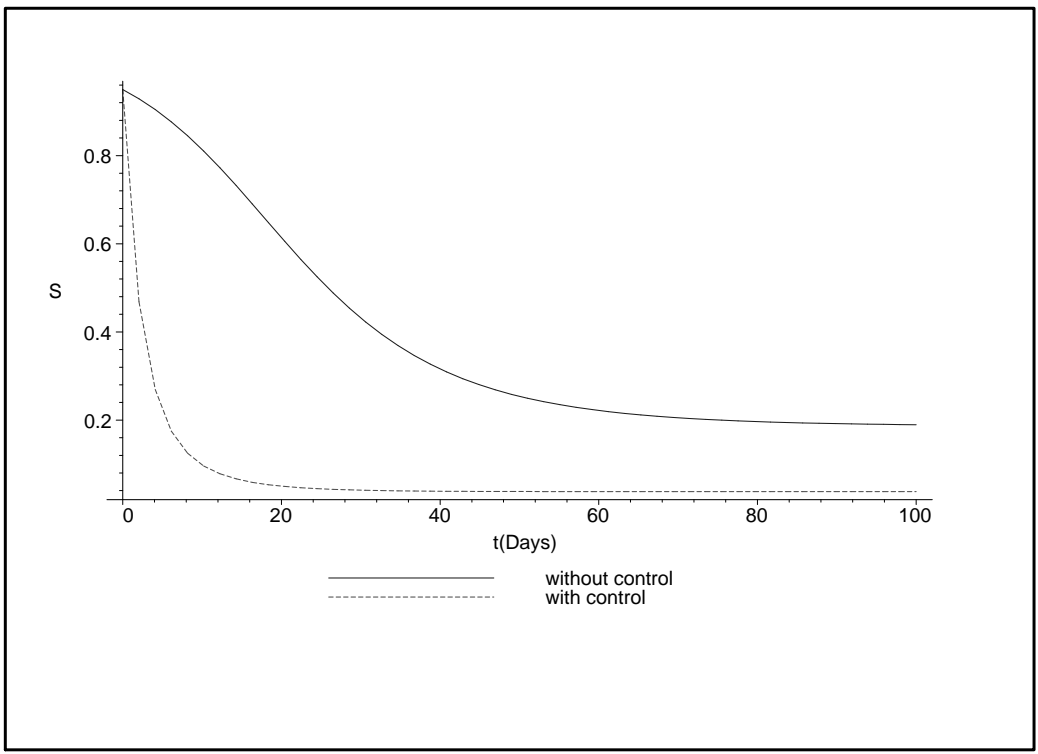}
\caption{Susceptible individuals $S(t)$
in case of optimal control \eqref{eq:E49} versus without control.}
\label{fig:7}
\end{figure}
\begin{figure}
\includegraphics[angle=-90,scale=0.6]{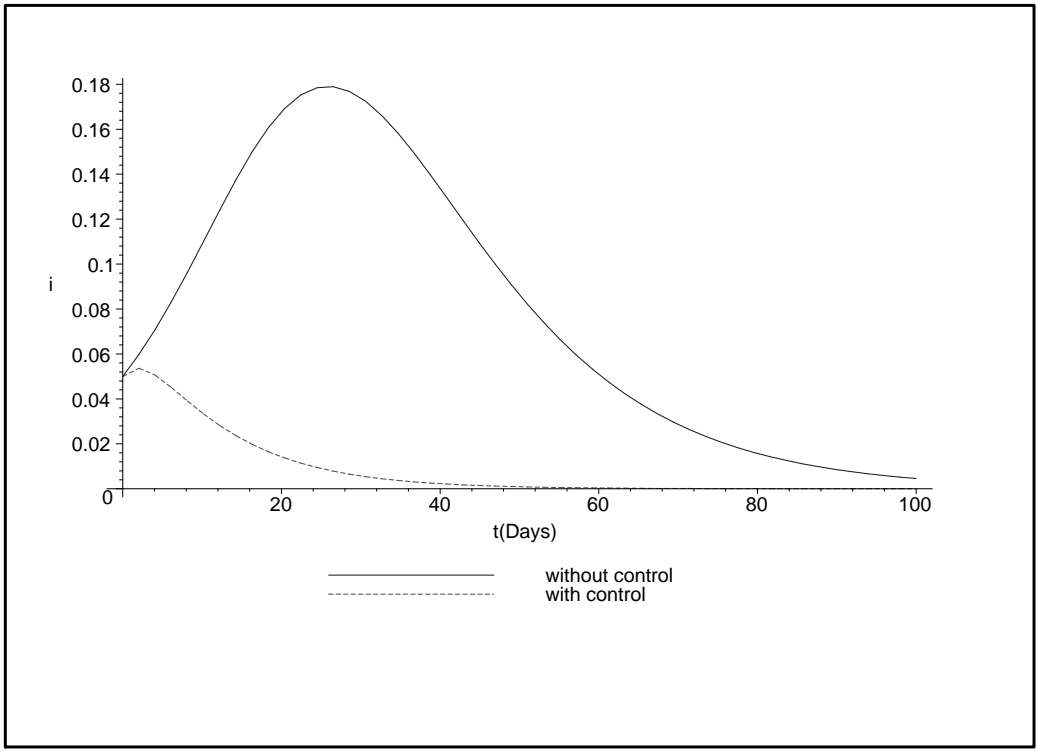}
\caption{Infected individuals $I(t)$
in case of optimal control \eqref{eq:E49}
versus without control.}
\label{fig:8}
\end{figure}
\begin{figure}
\includegraphics[angle=-90,scale=0.6]{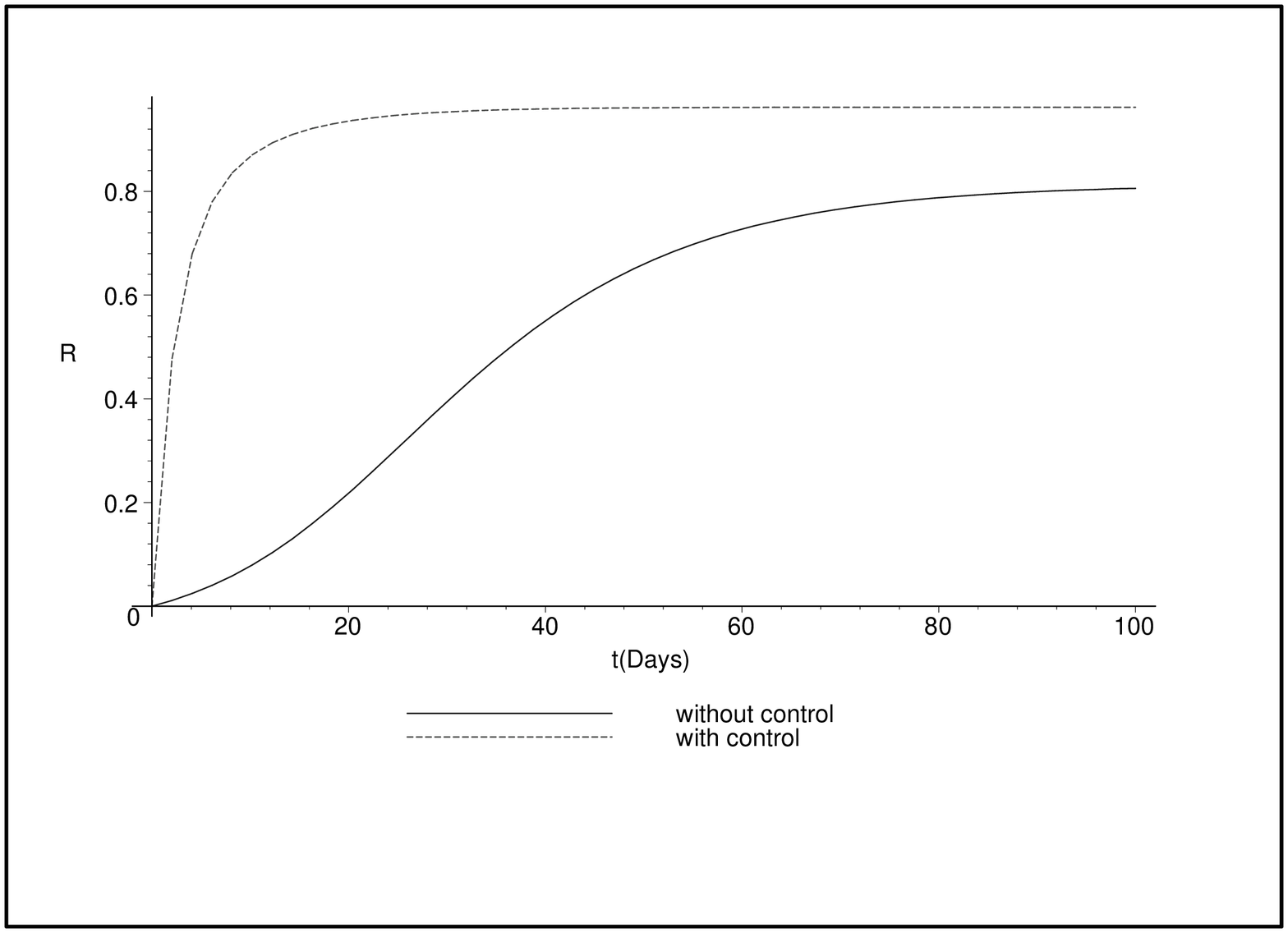}
\caption{Removed individuals $R(t)$ in case
of optimal control \eqref{eq:E49} versus without control.}
\label{fig:9}
\end{figure}


\section{Conclusions}
\label{conclu}

The analytical expression \eqref{eq:E49} for the optimal control
drawn at Figure~\ref{fig:6} is very similar to the optimal control
numerically depicted in \cite{MyID:321}. Similarly,
Figures~\ref{fig:7}, \ref{fig:8} and \ref{fig:9}, respectively
for the susceptible, infected and removed individuals, are
identical to the corresponding numerical results of \cite{MyID:321}.
We conclude that the numerical solutions found in \cite{MyID:321}
provide a good approximation to our analytical expressions.

We claim that the analytical method proposed here
can also be applied with success to other problems of optimal control
in mathematical epidemiology such as vector-borne, air-borne
and water-borne diseases. This question is under investigation
and will be addressed elsewhere.


\appendix


\section{Maple code}
\label{apen}

We have used the computer algebra system \textsf{Maple}
for all the computations. The reader interested in this
computer algebra system is referred, e.g., to \cite{MR2561903}.


\subsection{Maple code for Example~\ref{ex}}
\label{sec:ap:a1}

{\small
\begin{verbatim}
> restart:
> with(Physics):
> eq10:=J=Intc(diff(u(tau),tau)^(2)+1/2*A*u(tau)^2, tau);
> eq20:=Fundiff(eq10,u(t));
> eq20A:=dsolve({eq20,u(0)=U[0]});
> eq20B:=subs(_C2=U[0],eq20A);
> nas:=diff(S(t),t)=-beta*S(t)*i(t);
> nasB:=diff(i(t),t)=beta*S(t)*i(t)-u(t)*i(t);
> eq:=diff(i(t),t)=beta*S[0]*i(t)-rhs(eq20B)*i(t);
> eq1:=simplify(dsolve({eq,i(0)=iota[0]}),power,symbolic);
> eq2:=simplify(int(convert(series(rhs(eq1),t=0,2),polynom),t=0..T)
       + int((rhs(eq20B))^2*A/2,t=0..T),power,symbolic);
> eq3:=simplify(isolate(diff(eq2,U[0])=0,U[0]));
\end{verbatim}
}


\subsection{Maple code for the Ebola optimal control problem \eqref{eq:E4}--\eqref{eq:E6}}
\label{sec:ap:a2}

{\small
\begin{verbatim}
> restart:
> with(Physics):
> eq1:=J=Intc(diff(u(tau),tau)^(2)+1/2*A*u(tau)^2, tau);
> J(u) = Int([diff(u(tau),tau)^2+1/2*A*u(tau)^2],tau = 0 .. T);
> nis:=K(U[0])=Int(G(tau,A,U[0],Y[0]),tau=0..T)
         +(A/2)*int((U[0]*exp(-1/2*2^(1/2)*A^(1/2)*tau))^2,tau=0..T);
> nas:=diff(rhs(nis),U[0])=0;
> eq2:=Fundiff(eq1,u(t));
> eq4:=subs(_C2=U[0],dsolve({eq2,u(0)=U[0]}));
> restart:
> auxi:=u(t) = U[0]*exp(-1/2*2^(1/2)*A^(1/2)*t);
> aux0:=diff(s(t),t)=-rhs(auxi)*s(t);
> aux0A:=simplify(dsolve({aux0,s(0)=S[0]}),power,symbolic);
> aux0B:=s(t)=convert(series(rhs(aux0A),t=0,3),polynom);
> aux:=diff(i(t),t)=beta*s(t)*i(t)-mu*i(t);
> auxA:=subs(aux0A,aux);
> aux1:=dsolve({auxA,i(0)=iota[0]});
> aux1A:=i(t)=simplify(convert(series(rhs(aux1),t=0,5),polynom),power,symbolic);
> aux1B:=int(rhs(aux1A),t=0..T)+int(A*(rhs(auxi))^2/2,t=0..T);
> plas:=K(U[0])=subs(iota[0]=i[0],aux1B);
> plas1:=K(U[0])=i[0]*T+E[2]*T^2+E[3]*T^3
             +i[0]/48*E[4]*T^4-i[0]/240*E[5]*T^5
             +1/4*U[0]^2*2^(1/2)*A^(1/2)
             -1/4*2^(1/2)*A^(1/2)*U[0]^2*exp(-2^(1/2)*A^(1/2)*T);
> aux1C:=diff(aux1B,U[0])=0;
> aux1D:=isolate(aux1C,U[0]);
> yiyi:=U[0]*exp(-1/2*2^(1/2)*A^(1/2)*T/2)=U[0]/Q;
> isolate(U[0]*exp(-1/2*2^(1/2)*A^(1/2)*T/2)=U[0]/Q,A);
> param:={mu=0.1,beta=0.2,iota[0]=0.05,S[0]=0.95,T=100,Q=500};
> solu:=evalf(subs(param,aux1D));
> plot(rhs(solu),A=0.001..0.1);
> solu1:=evalf(isolate(U[0]*exp(-1/2*2^(1/2)*A^(1/2)*100/2)=U[0]/500,A));
> solu2:=evalf(subs(solu1,solu));
> plot(subs(solu2,solu1,rhs(auxi)),t=0..100);
> u:=subs(solu2,solu1,rhs(auxi));
> plot(u,t=0..100);
> with(plots):
> beta:=0.2;
> mu:=0.1;
> sysnc := diff(s(t),t)=-beta*s(t)*i(t),diff(i(t),t)=beta*s(t)*i(t)-mu*i(t),
                diff(r(t),t)=mu*i(t):
> fcns := {s(t),i(t),r(t)}:
> p:= dsolve({sysnc,s(0)=0.95,i(0)=0.05,r(0)=0},fcns,type=numeric,method=classical):
> odeplot(p, [[t,s(t)],[t,i(t)],[t,r(t)]],0..100);
> g:=odeplot(p, [[t,r(t)]],0..100,color=blue):
> gA:=odeplot(p, [[t,s(t)]],0..100,color=blue):
> gB:=odeplot(p, [[t,i(t)]],0..100,color=blue):
> sysc := diff(s(t),t)=-beta*s(t)*i(t)-u*s(t),
              diff(i(t),t)=beta*s(t)*i(t)-mu*i(t),
              diff(r(t),t)=mu*i(t)+u*s(t):
> fcns := {s(t),i(t),r(t)}:
> pc:= dsolve({sysc,s(0)=0.95,i(0)=0.05,r(0)=0},fcns,type=numeric,method=classical):
> sysc;
> odeplot(pc, [[t,s(t)],[t,i(t)],[t,r(t)]],0..50);
> g1:=odeplot(pc, [[t,r(t)]],0..100,color=red):
> g1A:=odeplot(pc, [[t,s(t)]],0..100,color=red):
> g1B:=odeplot(pc, [[t,i(t)]],0..100,color=red):
> display(g,g1);
> display(gA,g1A);
> display(gB,g1B);
\end{verbatim}
}


\section*{Acknowledgements}

This research was partially supported by the
Center for Research and Development in Mathematics and Applications (CIDMA)
within project UID/MAT/04106/2013 and by the Portuguese Foundation
for Science and Technology (FCT) through project TOCCATA, reference
PTDC/EEI-AUT/2933/2014. The authors are grateful to two anonymous
referees for valuable comments, suggestions and questions, which significantly
contributed to the quality of the paper.


\small



\end{document}